\documentclass[a4paper,12pt]{amsart}

\usepackage{amssymb,xspace,amscd,yhmath,mathdots,txfonts,
mathrsfs,calligra}

\usepackage[matrix,arrow,curve]{xy}\CompileMatrices

\usepackage{hyperref}

\usepackage{yfonts}[1998/10/03]

\DeclareMathAlphabet{\mathpzc}{OT1}{pzc}{m}{it}

\DeclareMathAlphabet{\mathcalligra}{T1}{calligra}{m}{n}
\DeclareFontShape{T1}{calligra}{m}{n}{<->s*[2.2]callig15}{}

\theoremstyle{plain}

\newtheorem{thm}{Theorem}

\newtheorem{lem}{Lemma}

\theoremstyle{definition}

\newtheorem{rmk}{Remark}

\DeclareMathAlphabet{\mathpzc}{OT1}{pzc}{m}{it}

\DeclareMathOperator{\R}{{\mathbb{R}}}

\DeclareMathOperator{\C}{\mathbb{C}}

\newcommand\sfB{\textsf{B}}

\newcommand\Lq{\mathpzc{L}}

\newcommand\sfX{\mathsf{X}}

\newcommand\Ci{\mathcal{C}^\infty}

\newcommand\Jd{\mathrm{J}}

\newcommand\pct{\mathrm{p}}

\newcommand\ttt{\texttt{t}}

\newcommand\sfS{\textsf{S}}

\newcommand\sfM{\textsf{M}}
\newcommand\sfT{\textsf{T}}
\newcommand\sfH{\textsf{H}}

   \def\DHLhksqrt#1#2{\setbox0=\hbox{$#1\sqrt{#2\,}$}\dimen0=\ht0
     \advance\dimen0-0.2\ht0
     \setbox2=\hbox{\vrule height\ht0 depth -\dimen0}
     {\box0\lower0.4pt\box2}}

\title
{Strictly pseudoconvex $CR$ structures on $\sfS^5$ that can be turned inside out}
\author[C.D.~Hill]{C. Denson Hill}
\address{C.D.\ Hill:
Department of Mathematics\\ Stony Brook University
\\ Stony Brook, N.Y. 11794 (USA)}
\email{dhill@math.stonybrook.edu}

\author[M.~Nacinovich]{Mauro Nacinovich}
\address{Mauro Nacinovich:
Dipartimento di Matematica\\ II Universit\`a di Roma
``Tor Ver\-ga\-ta''\\ Via della Ricerca Scientifica\\ 00133 Roma
(Italy)}
\email{nacinovi@mat.uniroma2.it}

\date{today}


\date{\today}
\begin{document}
\maketitle

\begin{abstract}
The existence of a complex structure on the real sphere $\sfS^{6}$
implies that of strictly pseudoconvex structures on $\sfS^{5}$
which can be turned inside out by a homotopy.
\end{abstract}

\section*{Introduction}

A perhaps locally defined, smooth, real codimension one, hypersurface $\sfM$ in a complex
manifold $\sfX$  inherits a $CR$ structure from the ambient complex structure on $\sfX$. If that $CR$
structure on $\sfM$ is strictly pseudoconvex, then $\sfM$ has two different sides that are distinguished
from one another by geometric considerations. One is the pseudoconvex side, and the other
is the pseudoconcave side. This distinction between the two sides is invariant with
respect to local biholomorphic diffeomorphisms of the ambient complex structure on $\sfX$,
and to local CR diffeomorphisms of M. This distinction between the two sides can be
understood in several different but equivalent ways:

Let x be a point on M. Then, as is well known, there is a local biholomorphic diffeomorphism
of a neighborhood of x in $\sfX$, which transforms $\sfM$ into an elementary strictly convex hypersurface.
The pseudoconvex side is then just the convex side, and the pseudoconcave side is just the
other concave side. That makes it easy to see that there is no local holomorphic disk, of complex
dimension equal to the CR dimension of M, whose center is at x, but otherwise stays on the
pseudoconvex side. On the other hand, there is a local holomorphic disk, of complex dimension
equal to the CR dimension of M, whose center is at x, but otherwise stays on the pseudoconcave
side of M. This can also be understood in terms of the scalar Levi form of $\sfM$ at x (see \cite{BHLN2020}), whose
signature is a local biholomorphic, and local CR diffeomorphism invariant.

For example, let $\sfM$ be the standard unit sphere $\sfS^{5}$, as it sits in  $\C^{3}$. Then, with the induced CR
structure coming from  $\C^{3}$, $\sfS^{5}$ is strictly pseudoconvex, the pseudoconvex side is the side of the
ball  $\sfB^{6}$ which the $\sfS^{5}$ bounds, and the pseudoconcave side is the exterior side. At every point
$x$ on $\sfS^{5}$, there is a local complex 2-dimensional disk whose center is at $x$, but otherwise stays on
the exterior side. There are no such disks on the interior side. This is the usual standard strictly
pseudoconvex  $CR$ structure one has on $\sfS^{5}$.

Now suppose that, by not changing the orientation of the $\sfS^{5}$ itself, one could smoothly deform
the  $CR$ structure on the fixed $\sfS^{5}$, through a one parameter family of  $CR$ structures on $\sfS^{5}$, so as
to switch the touching disks from the side on which they initially appear to the other side. That
would seem to be impossible for this standard  $CR$ structure on $\sfS^{5}$.

However, in this note we show that one can find strictly pseudoconvex  $CR$ structures on the fixed
$\sfS^{5}$, for which the touching disks do, in fact, switch sides in just that way. We refer to this by saying
that these  $CR$ structures on $\sfS^{5}$ can be "turned inside out". 
\par
We conjecture that for spheres  $\sfS^{k}$,
this phenomenon of having a strictly pseudoconvex  $CR$ stucture on  $\sfS^{k}$, that can be turned inside
out, occurs only when $k = 5$.

This also implies that there are  $CR$ structures on $\sfS^{5}$
with scalar Levi form degenerating at some point, all embeddable in a
complex $\sfS^{6}$.

\section*{A homotopy of $CR$ structures on $\sfS^{5}$}
Fix a complex structure  $\Jd$ 
on $\sfS^{6}$ (see \cite{Engel2026, Etesi2024}).
We identify 
$\sfS^{6}$ with the unit sphere: 
$$\sfS^{6}\,{\simeq}\,\left\{x\,{\in}\,\R^{7}\,\left|\,\|x\|^{2}{=}\,{\sum}_{h=0}^{6}x_{h}^{2}{=}\,1\right\}\right.$$
of $\R^{7}$
 and set 
 $$\sfM_{\ttt}\,{\coloneqq}\,\{x\,{\in}\,\sfS^{6}\,{\mid}\,x_{0}{=}\,\ttt\}.$$
For ${-}1\,{<}\,\ttt\,{<}\,1$, this is  a smooth real $5$-dimensional
hypersurface of $\sfS^{6}$ and  inherits a  $CR$ structure from $(\sfS^{6},\Jd)$. Let $\tauup$ be the restriction
of $dx_{0}$ to $\sfS^{6}$. 
If ${-}1\,{<}\,\ttt\,{<}\,1$ and $x\,{\in}\,\sfM_{\ttt}$, then $\tauup_{x}$ is a nonzero conormal to $\sfM_{\ttt}$
at $x$ and $\Jd^{*}\tauup_{x}$ generates the fiber  $\sfH^{0}_{x}\sfM_{\ttt}$ of the
characteristic bundle of $\sfM_{\ttt}$ at $x$. 
\par
We recall that, for a point $x$ of $\sfM_{\ttt}$,  
the \textsl{scalar} Levi form (see \cite{BHLN2020}) at $\xiup\,{\in}\,\sfH^{0}_{x}\sfM_{\ttt}$ is 
the Hermitian symmetric form 
\begin{gather*}
 \Lq_{\xiup}(Z_{1},Z_{2})=(2i)^{-1}\xiup([\tilde{Z}_{1},\bar{\tilde{Z}}_{2}]),\;\; \text{for}\; Z_{1},Z_{2}\in\sfT^{1,0}_{x}\sfM_{\ttt},
 \intertext{where $\tilde{Z}_{1},\tilde{Z}_{2}\in\Gamma^{\infty}(\sfM_{\ttt},\sfT^{1,0}\sfM_{\ttt})$ and $\tilde{Z}_{i}(x)=Z_{i}$,
 which, for $\xiup\,{=}\,\Jd^{*}\tauup_{x}$ yields}
 \Lq_{\Jd^{*}\tauup_{x}}(Z_{1},Z_{2})=Z_{1}\bar{Z}_{2}x_{0},\;\; \text{for}\; Z_{1},Z_{2}\in\sfT^{1,0}_{x}\sfM_{\ttt}.
 \end{gather*}
Take $x_{1},\hdots,x_{6}$ as \textit{real} smooth coordinates on $\sfS^{6}$ on neighbourhoods of the
south pole $\pct_{-}{=}({-}1,0,\hdots,0)$ and of the north pole $\pct_{+}{=}({+}1,0,\hdots,0)$. Then 
\begin{equation*}\begin{cases}
 x_{0}={-}\sqrt{1-x_{1}^{2}-\cdots-x_{6}^{2}}\simeq {-}1{+}\tfrac{1}{2}{\sum}_{h=1}^{6}x_{i}^{2} \;\;\text{near $\pct_{-}$},\\
 x_{0}={+}\sqrt{1-x_{1}^{2}-\cdots-x_{6}^{2}}\simeq {+}1{-}\tfrac{1}{2}{\sum}_{h=1}^{6}x_{i}^{2} \;\;\text{near $\pct_{+}$},
 \end{cases}
\end{equation*}
shows that the real valued function $x_{0}$ has positive definite Hessian at its critical point $\pct_{-}$
and negative definite Hessian at its critical point $\pct_{+}$. Taking holomorphic coordinates
in $(\sfS^{6},\Jd)$ near $\pct_{\pm}$, 
by the invariance 
of the signature of the Hessian at critical points, 
we get 
\begin{lem}
 The real valued function $x_{0}|_{\sfS^{6}}$ 
 is strictly pseudoconvex on a neighbourhood of $\pct_{-}$ and strictly pseudoconcave
 on a neighbourhood of~$\pct_{+}$. In particular, we can find $\epsilonup\,{>}\,0$ such that $\sfM_{\ttt}$
 is a smooth strictly pseudoconvex $CR$ hypersurface for real $\ttt$ with $0\,{<}\,(1\,{-}\,\ttt^{2})\,{<}\,\epsilonup$. \qed
\end{lem}
\begin{thm}
 $\sfS^{5}$
can be endowed with a strictly pseudoconvex $CR$ structure of hypersurface
type which can be "turned inside out", via a homotopy of $CR$ structures on $\sfS^5$, all
of which are of hypersurface type. 
\end{thm} 
\begin{proof} Let $y\,{=}\,(y_{1},\hdots,y_{6})$ be Euclidean coordinates on $\R^{6}$. For real $\thetaup$,
with $0\,{<}\,\thetaup\,{<}\,\piup$, we denote by $\sfT^{1,0}_{\thetaup}\sfS^{5}$ the $CR$ structure
on $\sfS^{5}$   
which is the pullback of the $CR$ structure 
$\sfT^{1,0}\sfM_{\cos(\thetaup)}\,{=}\, (\sfT^{1,0}\sfS^{6})\,{\cap}\,(\C\!\sfT\sfM_{\cos(\thetaup)})$
of the real hypersurface $\sfM_{\cos(\thetaup)}$ of $\sfS^{6}$
by   the parametrisation 
\begin{gather*} \phiup_{\thetaup}:
\sfS^{5}\ni{y}\to (\cos(\thetaup),\sin(\thetaup)\cdot{y})\in\sfM_{\cos(\thetaup)})\subset\sfS^{6}.
\intertext{For $0<\thetaup_{1}<\thetaup_{2}<\piup$,  
the map}
 [\thetaup_{1},\thetaup_{2}]\times\sfS^{5}\ni (\thetaup,y)\to (\cos(\thetaup),\sin(\thetaup)\cdot{y}) \in\sfS^{6},
\end{gather*}
yields  a homotopy between the $CR$ structures $\sfT^{1,0}_{\thetaup_{1}}\sfS^{5}$ and 
$\sfT^{1,0}_{\thetaup_{2}}\sfS^{5}$. We note that $\xiup_{\thetaup}{\coloneqq}\,\phiup_{\thetaup}^{*}(\Jd^{*}\tauup)$ 
is a nonzero one form of the characteristic bundle $\sfH^{0}_{\thetaup}\sfS^{5}$ of $(\sfS^{5},\sfT^{1,0}_{\thetaup}\sfS^{5})$.
If we take $\thetaup_{1},\thetaup_{2}$ with $0\,{<}\,\thetaup_{1}\,{<}\, \epsilonup$ and
$0\,{<}\,(\piup\,{-}\,\thetaup_{2})\,{<}\epsilonup$ for $0\,{<}\,\epsilonup\,{<\!<}1$, then by the previous lemma
the Levi form $\Lq_{\thetaup_{\thetaup}}$ on $(\sfS^{5},\sfT^{1,0}_{\thetaup}\sfS^{5})$
is negative definite at every point for $\thetaup\,{=}\,\thetaup_{1}$ and
positive definite at every point for $\thetaup\,{=}\,\thetaup_{2}$. 
\end{proof}

 The open set $\{x\,{\in}\,\sfS^{6}\,{\mid}\, x_{0}\,{<}\,\cos(\thetaup_{2})\}$ is a strictly pseudoconvex domain
 bounded by $\sfM_{\cos(\thetaup_{2})}$, while 
 $\{x\,{\in}\,\sfS^{6}\,{\mid}\, x_{0}\,{<}\,\cos(\thetaup_{1})\}$ is a strictly pseudoconcave domain
 bounded by $\sfM_{\cos(\thetaup_{1})}$.\par
 Thus, in going from $\thetaup_{1}$  to $\thetaup_{2}$, the touching complex 2-di\-men\-sional local
 holomorphic disks have changed sides; i.e., the $CR$ structure on $\sfS^{5}$ has been turned inside out.
  \begin{rmk}
By the construction, 
 these CR structures on $\sfS^5$  are all globally $CR$-embeddable in a complex $\sfS^6$.
\end{rmk}
\textbf{Conjecture}: Among spheres $\sfS^k$, of real dimension $k$, the only one that admits a
strictly pseudoconvex $CR$-hypersurface  structure, that can be turned inside out, is $\sfS^5$.
\par\smallskip
If $\ttt_{0}{=}\,\sup\{\ttt\,{\in}\,({-}1,1)\,{\mid}\,\Lq_{\xiup_{x}}\,{>}\,0,\; \forall x\in\sfM_{\ttt}\}$,
then ${-}1\,{<}\,\ttt_{0}\,{<}\,1$ 
and $\sfM_{\ttt_{0}}$ is weakly pseudoconvex $CR$ and its  Levi form degenerates at
some point.

\begin{rmk} In \cite{MV2002} it is shown that $\sfS^{5}$ admits a \textsl{Levi flat} $CR$ stucture.
On the other hand, a theorem of A. Haeflinger (see \cite{J1976}) implies that no
$\sfS^{5}$ may carry a Levi-flat \textit{real analytic} $CR$ structure. This in particular implies
that, if the embedding of the complex $\sfS^{6}$ in $\R^{7}$ is real analytic, then 
none of our $\sfM_{\ttt}$, for ${-}1\,{<}\,\ttt\,{<}\,1$, can be Levi flat. A more general question would be:
Can there be a $\Ci$ smooth  Levi flat real hypersurface,
which is not simply connected, in some complex $\sfS^6$?
\end{rmk}

\bibliographystyle{amsplain}
\renewcommand{\MR}[1]{}
\bibliography{cr}

\end{document}